\documentclass[12pt]{amsart}
\newtheorem{theorem}{Theorem}
\newtheorem{lemma}{Lemma}
\newtheorem{proposition}{Proposition}
\begin{document}
\title{Subelliptic biharmonic maps}
\author{}
\maketitle \centerline{{\large Sorin
Dragomir}\footnote{Universit\`a degli Studi della Basilicata,
Dipartimento di Matematica e Informatica, Contrada Macchia Romana,
85100 Potenza, Italy, e-mail: {\tt sorin.dragomir@unibas.it}} \and
{\large Stefano Montaldo}\footnote{Universit\`a degli Studi di
Cagliari, Dipartimento di Matematica e Informatica, Via Ospedale 72, 09124
Cagliari, Italy, e-mail: {\tt montaldo@unica.it}}}

\begin{abstract} We study subelliptic biharmonic maps i.e. smooth maps $\phi : M \to N$ from a compact strictly pseudoconvex CR manifold $M$
into a Riemannian manifold $N$ which are critical points of the energy functional $E_{2,b} (\phi ) = \frac{1}{2} \int_M \| \tau_b (\phi ) \|^2 \;\theta \wedge
(d \theta )^n$. We show that $\phi : M \to N$ is a
subelliptic biharmonic map if and only if its vertical lift $\phi \circ \pi : C(M) \to N$ to the (total space of the) canonical
circle bundle $S^1 \to C(M) \stackrel{\pi}{\longrightarrow} M$ is a biharmonic map with respect to the Fefferman metric $F_\theta$ on $C(M)$.
\end{abstract}

\section{Introduction}
Biharmonic maps were introduced by J. Eells \& L. Lemaire, \cite{EeLe}, as critical points $\phi \in C^\infty (M, N)$ of the bienergy functional
\[ E_2 (\phi ) = \frac{1}{2} \int_M \| \tau (\phi ) \|^2 \, d \, v_g  \]
where $M$ and $N$ are Riemannian manifolds and $\tau (\phi )$ is the tension field of $\phi : M \to N$. Biharmonic maps were
further investigated by G.~Jiang, \cite{JiGu}, who derived the first and second variation formulas for $E_2 (\phi )$
and gave several applications to the geometry of the second fundamental form of a submanifold in a Riemannian manifold.
In the last decade a large amount of work has been
devoted to biharmonic maps with particular attention to the constructions and classifications of proper biharmonic maps and proper biharmonic submanifolds, 
see e.g. \cite{BaKa}-\cite{BaFaOu}, \cite{BMO08}-\cite{BMO10}, \cite{CaMoOn1}-\cite{CaMoOn2}, \cite{IcInUr}, \cite{O10}-\cite{O-pre}.
\par
A program aiming to extending results on nonlinear elliptic systems of variational origin to the hypoelliptic case was started by J. Jost \& C-J. Xu, \cite{JoXu}.
Given a H\"ormander system of vector fields $\{ X_1 , \cdots , X_p \}$ on an open set $U \subset {\mathbb R}^m$ and a map $\phi \in C^\infty (U, N)$ the function
\[ e(\phi ) = \frac{1}{2} \sum_{a=1}^p X_a (\phi^j ) X_a (\phi^k ) (h_{jk} \circ \phi ) \]
is globally defined and generalizes the ordinary energy density of $\phi$. Here $h_{jk}$ is the Riemannian metric on $N$ in a local coordinate system
$(V, y^j )$ and $\phi^j = y^j \circ \phi$. A subelliptic harmonic map is a critical point $\phi \in C^\infty (U, N)$ of the functional (cf.
\cite{JoXu})
\[ E_X (\phi ) = \int_\Omega e(\phi ) \; d x \]
where $\Omega \subset {\mathbb R}^m$ is a bounded domain such that $\overline{\Omega} \subset U$. The Euler-Lagrange equations of the variational
principle $\delta \, E_X (\phi ) = 0$
are
\[ - H \phi^i + \left( \Gamma^i_{jk} \circ \phi \right) \sum_{a=1}^p X_a (\phi^j ) X_a (\phi^k ) = 0 \]
where $H \equiv \sum_{a=1}^p X_a^\ast X_a$ is the H\"ormander operator and
$\Gamma^i_{jk}$ are the Christoffel symbols of $h_{jk}$.
Subelliptic harmonic maps were recognized (cf.
E. Barletta \& S. Dragomir \& H. Urakawa, \cite{BaDrUr1}) as the local manifestation
of pseudoharmonic maps i.e. critical points $\phi \in C^\infty (M, N)$ of the functional
\[ E_{1,b}(\phi ) = \frac{1}{2} \int_M {\rm trace}_{G_\theta} \left( \Pi_H \, \phi^\ast h \right) \; \theta \wedge (d \theta )^n \, . \]
Here $M$ is a compact strictly pseudoconvex CR manifold, of CR dimension $n$, and $\theta$ is a contact form on $M$ such that the Levi form $G_\theta$ is positive
definite. Also $h$ is the Riemannian metric on $N$ and $\Pi_H \, \phi^\ast h$ is the restriction of the bilinear form $\phi^\ast h$ to the Levi, or maximally
complex, distribution $H(M)$. The Euler-Lagrange equations of $\delta \, E_{1, b}(\phi ) = 0$ may be written as
$\tau_b (\phi ) = 0$ where the field $\tau_b (\phi ) \in C^\infty (\phi^{-1} T(N))$ is locally given by
\[ \tau_b (\phi )^i = \Delta_b \phi^i + \sum_{a=1}^{2n} X_a (\phi^j ) X_a (\phi^k ) \left( \Gamma^i_{jk} \circ \phi \right) . \]
Here $\{ X_a : 1 \leq a \leq 2n \}$ is a local $G_\theta$-orthonormal frame
of $H(M)$. Also $\Delta_b$ is the sublaplacian i.e. the second order differential operator given by
\[ \Delta_b u = {\rm div} \left( \nabla^H u \right) , \;\;\; u \in C^\infty (M). \]
The divergence operator is meant with respect to the volume form $\theta \wedge (d \theta )^n$ while the horizontal gradient
is given by $\nabla^H u = \Pi_H \nabla u$ and $g_\theta (\nabla u , X) = X(u)$ for any $X \in \mathfrak{X}(M)$. CR manifolds occur mainly as boundaries
$M = \partial \Omega$ of domains $\Omega \subset {\mathbb C}^{n+1}$ and boundary values of Bergman-harmonic maps $\Phi : \Omega \to N$ may be shown (cf.
\cite{DrKa}) to be pseudoharmonic provided $\Phi$ has vanishing normal derivatives (thus motivating our use of the index $b$ for "boundary" analogs to geometric
objects such as the tension field $\tau (\Phi )$, the second fundamental form $\beta (\Phi )$, etc.). The similar boundary behavior of biharmonic maps from a
strictly pseudoconvex domain
$\Omega \subset {\mathbb C}^{n+1}$ endowed with the Bergman metric is unknown.
\par
The approaches in \cite{JoXu} and \cite{BaDrUr1}
overlap partially, as follows. For any local $G_\theta$-orthonormal frame $\{ X_a : 1 \leq a \leq 2 n \}$ in $H(M)$ defined on a local coordinate neighborhood
$\varphi : U \to {\mathbb R}^{2n+1}$ the push-forward $\{ (d \varphi ) X_a : 1 \leq a \leq 2 n \}$ is a H\"ormander system on $\varphi (U)$ and $\Delta_b = - H$.
However in \cite{JoXu} and in general in the theory of H\"ormander vector fields
the Euclidean dimension $m$ is arbitrary (as opposed to $m = 2n+1$ in the CR case), the vector fields forming the given system are allowed to be linearly dependent
(at particular points), and the formal adjoints $X_a^\ast$ are meant with respect to the Euclidean metric on $\Omega$ (rather than the Webster metric - a non flat
Riemannian metric springing from the given CR structure in the presence of a contact form). It is therefore a natural problem,
within J.Jost \& C-J. Xu's program, to study critical points $\phi \in C^\infty (M, N)$ of the functional
\[ E_{2, b} (\phi ) = \frac{1}{2} \int_M \| \tau_b (\phi ) \|^2 \, \theta \wedge (d \theta )^n \, . \]
These are referred to as {\em subelliptic biharmonic maps}. We give a geometric interpretation of subelliptic biharmonic maps in terms of ordinary biharmonic
maps from a Lorentzian manifold (the total space of the canonical circle bundle $S^1 \to C(M) \to M$ endowed with the Fefferman metric $F_\theta$, \cite{DrTo}).
The paper is organized as follows. In $\S \, 2$ we consider ordinary biharmonic maps from a Fefferman space-time. The rough Laplacian (a degenerate
elliptic operator appearing in the principal part
of the subelliptic biharmonic map system) is discussed in $\S \, 3$. The first variation formula for the functional $E_{2,b} : C^\infty (M, N) \to {\mathbb R}$
is derived in $\S \, 4$ while the main result (cf. Theorem \ref{t:SubBih1} below) is proved in $\S \, 5$. A few open problems are outlined in $\S \, 6$.

\section{Biharmonic maps from Fefferman space-times}
Let $(M, T_{1,0}(M))$ be a compact orientable CR manifold, of CR dimension $n$, where $T_{1,0}(M)$ is its CR
structure. Let us assume that $M$ is strictly pseudoconvex and $\theta$ a contact form on $M$ such that the Levi form
$G_\theta (X,Y) = (d \theta )(X, J Y)$, $X,Y \in H(M)$, is positive definite. Here $H(M) = {\rm Re} \left\{ T_{1,0}(M)
\oplus T_{0,1}(M) \right\}$ is the Levi distribution of $M$ and $J : H(M) \to H(M)$, $J(Z +
\overline{Z}) = i(Z - \overline{Z})$, $Z \in T_{1,0}(M)$, its complex structure ($i = \sqrt{-1}$).
For all needed notions of CR and pseudohermitian geometry we rely on \cite{DrTo}. Let $T$ be the Reeb vector of $(M, \theta )$ i.e. the
nowhere zero globally defined tangent vector field field transverse to $H(M)$ determined by $\theta (T) = 1$ and $T \,
\rfloor \, d \theta = 0$. Let $\nabla$ be the Tanaka-Webster connection of $(M, \theta )$ (cf. Theorem
3.1 in \cite{DrTo}, p. 25, for the axiomatic description of $\nabla$).
\par
Let $S^1 \to C(M) \stackrel{\pi}{\longrightarrow} M$ be the canonical circle bundle over $M$ and $F_\theta$ the Fefferman
metric on $C(M)$ (cf. J.M. Lee, \cite{Lee}, or Definition 2.15 in \cite{DrTo}, p. 128) i.e. the Lorentzian metric on $C(M)$ given by
(cf. (2.30) in \cite{DrTo}, p. 127)
\begin{equation} F_\theta = \pi^* \tilde{G}_\theta + 2\left( \pi^* \theta \right) \odot \sigma ,
\label{e:SubBih1}
\end{equation}
\begin{equation}
\sigma = \frac{1}{n+2} \left\{ d \gamma + \pi^* \left( i \, \omega^\alpha_\alpha - \frac{i}{2} \, g^{\alpha\overline{\beta}} d
g_{\alpha\overline{\beta}} - \frac{\rho}{4(n+1)} \, \theta \right) \right\} .
\label{e:SubBih2}
\end{equation}
As to the notation in (\ref{e:SubBih1})-(\ref{e:SubBih2}), given a local frame $\{ T_\alpha : 1 \leq \alpha \leq n \}$ of
$T_{1,0}(M)$ we set $g_{\alpha\overline{\beta}} = G_\theta (T_{\alpha} , T_{\overline{\beta}})$ (with $T_{\overline{\alpha}}
= \overline{T}_\alpha$). Also $\omega^\alpha_\beta$ are the connection $1$-forms of $\nabla$ with respect to $\{ T_\alpha : 1
\leq \alpha \leq n \}$ i.e. $\nabla T_\beta = \omega^\alpha_\beta \otimes T_\alpha$. Moreover $\rho = g^{\alpha\overline{\beta}}
R_{\alpha\overline{\beta}}$ is the pseudohermitian scalar curvature of $(M , \theta )$ (cf. e.g. \cite{DrTo}, p. 50). The
$(0,2)$-tensor field $\tilde{G}_\theta$ in (\ref{e:SubBih1}) is obtained by extending the Levi form $G_\theta$ to a degenerate
form defined on the whole of $T(M)$. By definition $\tilde{G}_\theta = G_\theta$ on $H(M) \otimes H(M)$ and
$\tilde{G}_\theta (X, T) = 0$ for any $X \in T(M)$. By a result of C.R. Graham, \cite{Gra}, the (globally defined) $1$-form $\sigma$
is a connection form in the canonical circle bundle. Let $X^\uparrow \in \mathfrak{X}(C(M))$ be the horizontal lift of $X \in \mathfrak{X}(M)$
with respect to $\sigma$. If $S \in \mathfrak{X}(C(M))$ is the tangent to the $S^1$ action then $T^\uparrow - S$ is a nowhere vanishing
globally defined timelike vector field hence $(C(M), F_\theta )$ is a time-oriented Lorentzian manifold, referred to as the {\em Fefferman space-time}.
$F_\theta$ was discovered by C. Fefferman, \cite{Fef}, in connection with the study of the boundary behavior of the Bergman kernel of a bounded
domain $\Omega \subset {\mathbb C}^n$ (as a Lorentzian metric on $C(\partial \Omega ) \approx \partial \Omega \times S^1$). An array of popular nonlinear
problems arise from $F_\theta$ e.g. the CR Yamabe problem (cf. \cite{DrTo}, p. 159-160) is the projection via $\pi : C(M) \to M$ of the ordinary
Yamabe problem for $F_\theta$. While the principal part of the Yamabe equation on $C(M)$ is the wave operator $\square$ (hence the Yamabe equation
on $C(M)$ is not elliptic) the principal part of the projected equation is the sublaplacian $\Delta_b$ (hence subelliptic theory applies, cf.
\cite{DrTo}, p. 176-210). Also pseudoharmonic maps from strictly pseudoconvex CR manifolds may be characterized as base maps of $S^1$-invariant
harmonic maps from $(C(M), F_\theta )$, thus suggesting that Theorem \ref{t:SubBih1} below should hold.
\par
Several basic facts in harmonic map theory are known to extend in a straightforward manner from the Riemannian to the semi-Riemannian setting (cf.
\cite{BaWo}, p. 427-452) e.g. a smooth map of semi-Riemannian manifolds has a well defined tension tensor field. In particular the applications we seek for are to maps
from Lorentzian to Riemannian manifolds. A $C^\infty$ map $\Phi : C(M) \to N$ into a real $\nu$-dimensional
Riemannian manifold $(N, h)$ is {\em biharmonic} if $\Phi$ is a critical point of the functional
\begin{equation}\label{e:SubBih3}
{\mathbb E}_2 (\Phi ) = \frac{1}{2} \int_{C(M)} \| \tau( \Phi )
\|^2 \; d \, {\rm vol}(F_\theta ).
\end{equation}
The {\em tension field} $\tau (\Phi )$ is the $C^\infty$ cross-section in the pullback bundle $\Phi^{-1} T N \to C(M)$ locally given by
\[ \tau (\Phi ) = \left( \square \Phi^i + \left( \Gamma^i_{jk} \circ \Phi \right) \frac{\partial \Phi^j}{\partial u^p}
\frac{\partial \Phi^k}{\partial u^q} F^{pq} \right) X^\Phi_i \]
where $\square$ is the Laplace-Beltrami operator of $F_\theta$ (the {\em wave operator} as $F_\theta$ is Lorentzian) and $\Phi^i
= y^i \circ \Phi$.  Also if $(U, x^A )$ is a local coordinate system on $M$ such that $\Phi^{-1}(V) \subset \pi^{-1} (U)$
and $\gamma : \pi^{-1} (U) \to {\mathbb R}$ is a
local fibre coordinate on $C(M)$ then $(\pi^{-1} (U), \; u^A = x^A \circ \pi , \; u^{2n+2} = \gamma )$ are the naturally induced
local coordinates on $C(M)$ and $[F^{pq}] = [F_{pq}]^{-1}$ while $F_{pq} = F_\theta (\partial_p , \partial_q )$. Here $\partial_p$
is short for $\partial /\partial u^p$. Finally $X_i^\Phi$ is the {\em natural lift} of $\partial_i = \partial /\partial y^i$ i.e.
the local smooth section in $\Phi^{-1} TN \to C(M)$ given by $X_i^\Phi (z) = \left( \partial_i \right)_{\Phi (z)}$ for any $z
\in \Phi^{-1}(V)$. The bundle metric $h^\Phi$ appearing in (\ref{e:SubBih3}) is naturally induced by $h$ in $\Phi^{-1} T N
\to C(M)$ so that $h^\Phi (X^\Phi_i , X^\Phi_j ) = h_{ij} \circ \Phi$.
Our main result is
\begin{theorem}
Let $M$ be a compact strictly pseudoconvex CR manifold, of CR dimension $n$, and $\theta$ a contact form on $M$ with $G_\theta$
positive definite. Let $\Phi : C(M) \to N$ be a smooth $S^1$-invariant map and $\phi : M \to N$ the corresponding base map.
Then ${\mathbb E}_2 (\Phi ) = 2 \pi E_{2,b} (\phi )$. Consequently if $\Phi$ is biharmonic then $\phi$
is a critical point of $E_2$. The Euler-Lagrange equations of the variational principle $\delta \, E_{2,b} (\phi ) = 0$ are
\begin{equation}
BH_b (\phi ) \equiv
\Delta_b^\phi \, \tau_b (\phi ) + {\rm trace}_{G_\theta} \left\{ \Pi_H \, \hat{R}^h (\tau_b (\phi ), \, \phi_* \, \cdot ) \phi_* \, \cdot \right\} = 0
\label{e:SubBih4}
\end{equation}
where $\Delta^\phi_b$ is the rough sublaplacian and $R^h$ the curvature tensor field of $N$. Consequently the vertical lift to
$C(M)$ of any $C^\infty$ solution $f$ to {\rm (\ref{e:SubBih4})} is a biharmonic map {\rm (}with respect to the Fefferman metric
$F_\theta${\rm )}.
\label{t:SubBih1}
\end{theorem}
Here $\tau_b ( \phi ) = \tau (\phi  ; \theta , \nabla^h )$ is the {\em subelliptic tension field} of $\phi : M \to N$ i.e. the $C^\infty$
section in $\phi^{-1} T N \to M$ locally given by
\[ \tau_b (\phi ) = \left\{ \Delta_b \phi^i + 2 g^{\alpha\overline{\beta}}
\left( \Gamma^i_{jk} \circ \phi \right) T_\alpha (\phi^j ) T_{\overline{\beta}} (\phi^k ) \right\} X_i^\phi \]
where $\Delta_b$ is the sublaplacian of $(M , \theta )$ (cf. the Introduction or Definition 2.1 in \cite{DrTo}, p. 134), $\phi^i = y^i \circ \phi$
and $X_i^\phi (x) = \left( \partial_i \right)_{\phi (x)}$ for any $x \in \phi^{-1}(V)$ (the natural lift of $\partial /\partial y^i$ as
a section in $\phi^{-1} T N \to M$). Note that $X_i^\Phi = X_i^\phi \circ \pi$. The {\em rough sublaplacian} is the second order
differential operator locally given by
\[ \Delta_b^\phi \, V = \sum_{a=1}^{2n} \{ (\phi^{-1} \nabla^h )_{X_a} (\phi^{-1} \nabla^h )_{X_a} V - (\phi^{-1} \nabla^h )_{\nabla_{X_a} X_a} V \} \]
for any $V \in C^\infty (\phi^{-1} T N)$, where $\{ X_a : 1 \leq a \leq 2n \}$ is a local orthonormal (i.e. $G_\theta (X_a ,
X_b ) = \delta_{ab}$) frame of $H(M)$. Moreover $\phi^{-1} \nabla^h$ is the connection in $\phi^{-1} T N \to M$ defined as the pullback of
the Levi-Civita connection $\nabla^h$ of $(N, h)$ by $\phi$ i.e.
\[ (\phi^{-1} \nabla^h )_{\partial /\partial x^A} X^\phi_k = \frac{\partial \phi^j}{\partial x^A} \left( \Gamma^i_{jk} \circ \phi \right) X^\phi_k . \]
Given a bilinear form $B$ on $T(M)$ we denote by $\Pi_H B$ the restriction of $B$ to $H(M) \otimes H(M)$. We shall need the following
\begin{theorem}
For any smooth map $\phi : M \to N$ the rough sublaplacian $\Delta^\phi_b$ is a formally self adjoint {\rm (}with respect to the $L^2$ inner
product $(V,W) = \int_M h^\phi (V,W) \, \theta \wedge (d \theta )^n$, $V,W \in C^\infty (\phi^{-1} T N)${\rm )} second order
differential operator locally expressed as
\begin{equation}
\Delta_b^\phi \, V = \left\{ \Delta_b V^i + 2 \sum_{a=1}^{2n} X_a (\phi^j ) \left( \Gamma^i_{jk} \circ \phi \right) X_a (V^k ) + \right.
\label{e:SubBih5}
\end{equation}
\[ \left. + \left[ \left( \Gamma^i_{jk} \circ \phi \right) \Delta_b \phi^j + \sum_{a=1}^{2n} X_a (\phi^j ) X_a (f^\ell )
\left( \frac{\partial \Gamma^i_{jk}}{\partial y^\ell} + \Gamma^s_{k\ell} \Gamma^i_{js} \right) \circ \phi \right] V^k \right\} X_i^\phi \]
where $V = V^i X_i^\phi$, for any local orthonormal frame $\{ X_a : 1 \leq a \leq 2 n \}$ of $H(M)$
defined on the open set $U \subset M$ and any local coordinate system $(V, y^i )$ on $N$ such that $\phi^{-1} (V) \subset U$. Let
$D^*$ be the formal adjoint of $D = (\phi^{-1} \nabla^h )^H$ i.e. $(D^* \varphi , V) = (\varphi , D V)$ for any $\varphi \in
C^\infty (H(M)^* \otimes \phi^{-1} T N)$ and any $V \in C^\infty_0 (\phi^{-1} T N)$. Then
\begin{equation}
\Delta_b^\phi = - D^* D.
\label{e:RoughSubLaplacian}
\end{equation}
In particular $(\Delta_b^\phi V , V) \leq 0$.
\label{p:SubBih1}
\end{theorem}
The proof of Theorem \ref{p:SubBih1} will be given in $\S \, 3$. If $N = {\mathbb R}^\nu$ then (by (\ref{e:SubBih4})-(\ref{e:SubBih5})) the subelliptic biharmonic
map equations become $L \phi^i = 0$ where $L \equiv \Delta_b \circ \Delta_b$ (the {\em bi-sublaplacian}) is a fourth order hypoelliptic operator.
The analysis of the scalar case $N =
{\mathbb R}$ (maximum principles, existence of Green functions for $L$, Harnack inequalities for positive solutions to $L u = 0$,
etc.) is however open. The calculation of the Green function for $\Delta^2 \equiv \Delta \circ \Delta$ (where $\Delta$ is the ordinary Laplacian on ${\mathbb
R}^n$) is due to T. Boggio, \cite{Bog}. The existence of the Green function for $\Delta_b$ follows from work by J.M. Bony,
\cite{Bon}, while estimates (on the Green function and its derivatives) were got by A. S\'anchez-Calle, \cite{San}, yet the
problem of adapting their techniques to the bi-sublaplacian is unsolved.

\section{The rough sublaplacian}
The second order differential operator $\Delta_b^\phi$ is similar to the rough Laplacian on vector fields due to G. Wiegmink,
\cite{Wie}, and C.M. Wood, \cite{Woo}. Let $\phi : M \to N$ be a smooth map and $\phi^{-1} T N \to M$ the pullback bundle. Let $h^\phi$
(respectively $\phi^{-1} \nabla^h$) be the pullback of the Riemannian metric $h$ (respectively of the Levi-Civita connection $\nabla^h$)
by $\phi$. Then $h^\phi$ is parallel with respect to $\phi^{-1} \nabla^h$.
We shall establish
\begin{lemma} For any $V, W \in C^\infty (\phi^{-1} T N)$ there is a smooth tangent vector field $X_\phi$ on $M$ such that
\begin{equation}
h^\phi (\Delta_b^\phi \, V  , \, W) = \Delta_b \left[ h^\phi (V,W) \right] + h^\phi (V , \, \Delta^\phi_b \, W) - 2 \, {\rm div}(X_\phi )
\label{e:SubBih7}
\end{equation}
where the divergence is taken with respect to the volume form $\Psi = \theta \wedge (d \theta )^n$ i.e. ${\mathcal L}_{X_\phi} \Psi
= {\rm div}(X_\phi ) \Psi$.
\label{l:SubBih1}
\end{lemma}
{\em Proof}. Let $\{ X_a : 1 \leq a \leq 2 n \}$ be a local orthonormal frame of $H(M)$. As $(\phi^{-1} \nabla^h ) h^\phi = 0$
\[ h^\phi \left( (\phi^{-1} \nabla^h )_{X_a} (\phi^{-1} \nabla^h )_{X_a} V \, , \, W \right) = \]
\[ = X_a \left( h^\phi \left( (\phi^{-1} \nabla^h )_{X_a} V  \, , \, W \right) \right) - h^\phi \left( (\phi^{-1} \nabla^h )_{X_a} V \, , \,
(\phi^{-1} \nabla^h )_{X_a} W \right) = \]
\[ = X_a^2 \left( h^\phi \left( V,W \right) \right) - 2 X_a \left( h^\phi \left( V , \, (\phi^{-1} \nabla^h )_{X_a} W \right) \right) + h^\phi
\left( V, \, (\phi^{-1} \nabla^h )_{X_a}^2 W \right) , \]
\[ h^\phi \left( (\phi^{-1} \nabla^h )_{\nabla_{X_a} X_a} V \, , \, W\right) = \]
\[ = (\nabla_{X_a} X_a )(h^\phi (V,W)) - h^\phi (V, \, (\phi^{-1} \nabla^h )_{\nabla_{X_a} X_a} W ). \]
Let us recall that $\Delta_b u = {\rm div} \left( \nabla^H u \right)$, $u \in C^2 (M)$, where $\nabla^H u \in C^\infty (H(M))$ (the {\em horizontal gradient}
of $u$) is given by $\nabla^H u = \Pi_H \nabla u$ and $g_\theta (\nabla u , X) = X(u)$ for any $X \in \mathfrak{X}(M)$.
Also $\Pi_H : T(M) \to H(M)$ is the natural projection associated to the direct
sum decomposition $T(M) = H(M) \oplus {\mathbb R} T$ while $g_\theta$ is the {\em Webster metric} of $(M , \theta )$ (cf.
Definition 1.10 in \cite{DrTo}, p. 9). A large amount of the existing subelliptic theory is built on the Heisenberg group (cf. \cite{BoLaUg},
p. 155) i.e. on the noncommutative Lie group ${\mathbb H}_n \equiv {\mathbb C}^n \times {\mathbb R}$ with the multiplication law
\[ (z, t) \cdot (w , s) = (z+w, \; t + s + 2 \, {\rm Im}(z \cdot \overline{w})) \]
for any $(z,t), \, (w,s) \in {\mathbb H}_n$. Let $\overline{Z}_\alpha$ be the {\em Lewy operators} i.e.
\[ Z_\alpha \equiv \frac{\partial}{\partial z^\alpha} + i \overline{z}^\alpha \, \frac{\partial}{\partial t} \, , \;\;\; 1 \leq \alpha \leq n. \]
Then $T_{1,0}({\mathbb H}_n ) = \sum_{\alpha =1}^n {\mathbb C} Z_\alpha$ is a CR structure on ${\mathbb H}_n$, of CR dimension $n$, making ${\mathbb H}_n$
into a CR manifold (cf. e.g. \cite{DrTo}). The Levi distribution $H({\mathbb H}_n )$ is spanned by the left invariant vector fields $X_a \in \mathfrak{X}({\mathbb H}_n )$
given by
\[ X_\alpha \equiv \frac{\partial}{\partial x^\alpha} + 2 y^\alpha \frac{\partial}{\partial t} , \;\;\; X_{n + \alpha} \equiv \frac{\partial}{\partial y^\alpha}
- 2 x^\alpha \frac{\partial}{\partial t} , \]
and the horizontal gradient is familiar (cf. \cite{BoLaUg}, p. 68) in subelliptic theory as the horizontal $H$-gradient i.e. $\nabla^H u =
\sum_{a=1}^{2n} X_a (u) X_a$. Note also that for the Heisenberg group
$H = - \sum_{a=1}^{2n} X_a^2$ (H\"ormander's sum of squares of vector fields) as $X_a^\ast = - X_a$.
\par
Let $\nabla$ be the Tanaka-Webster connection of $(M, \theta )$. As $\nabla \Psi = 0$ the divergence of a vector field may also be computed as the trace
of its covariant derivative with respect to $\nabla$. Hence $\Delta_b u$ may be locally written as
\[ \Delta_b u  = \sum_{a=1}^{2n} \left\{ X_a^2 u - \left( \nabla_{X_a} X_a \right) u \right\} . \]
Consequently
\begin{equation}
h^\phi \left( \Delta_b^\phi \, V , \, W \right) = \Delta_b \left[ h^\phi (V,W) \right] + h^\phi \left( V, \, \Delta_b^\phi \, W \right) +
\label{e:SubBih8}
\end{equation}
\[ + 2 \sum_a \{ h^\phi (V, \, (\phi^{-1} \nabla^h )_{\nabla_{X_a} X_a} W) - X_a (h^\phi (V, \, (\phi^{-1} \nabla^h )_{X_a} W)) \} . \]
Let $X_\phi \in H(M)$ be the vector field determined by
\[ G_\theta (X_\phi , Y) = h^\phi (V, \, (\phi^{-1} \nabla^h )_Y W) \]
for any $Y \in H(M)$. Then (by $\nabla g_\theta = 0$)
\[ \sum_a X_a (h^\phi (V, \, (\phi^{-1} \nabla^h )_{X_a} W)) = \sum_a X_a (G_\theta (X_\phi , X_a )) = \]
\[ = \sum_a \{ g_\theta (\nabla_{X_a} X_f \, , \, X_a ) + g_\theta (X_\phi \, , \, \nabla_{X_a} X_a ) \} = \]
\[ = {\rm div}(X_\phi ) + \sum_a h^\phi (V, \, (\phi^{-1} \nabla^h )_{\nabla_{X_a} X_a} W) . \]
Indeed $g_\theta (\nabla_T X_\phi \, , \, T) = 0$ (as $H(M)$ is parallel with respect to $\nabla$).
Together with (\ref{e:SubBih8}) this leads to (\ref{e:SubBih7}). Lemma \ref{l:SubBih1} is proved.
\vskip 0.1in
Let us assume that either $M$ is compact or at least one of the sections $V,W$ has compact support. At this point we may integrate
(\ref{e:SubBih7}) over $M$ and use Green's lemma to show that $(\Delta_b^\phi \, V , \, W) = (V, \, \Delta_b^\phi \, W)$. Moreover, if
$V = V^i X_i^\phi$ is a $C^\infty$ section in $\phi^{-1} T N \to M$ then
\begin{equation}
(\phi^{-1} \nabla^h )_X V = \{ X(V^i ) + X(\phi^j ) V^k \left( \Gamma^i_{jk} \circ \phi \right) \} X_i^\phi
\end{equation}
for any $X \in {\mathcal X}(M)$. The proof of (\ref{e:SubBih5}) follows from
\[ (\phi^{-1} \nabla^h )_X^2 V = \left\{ X^2 (V^i ) + 2 X(\phi^j ) \left( \Gamma^i_{jk} \circ f \right) X(V^k ) + \right. \]
\[ \left. + \left[ X^2 (\phi^j ) \left( \Gamma^i_{jk} \circ \phi \right) + X(\phi^j ) X(f^\ell ) \left( \frac{\partial \Gamma^i_{jk}}{\partial y^\ell} +
\Gamma^m_{k\ell} \Gamma^i_{jm} \right) \circ \phi \right] V^k \right\} X_i^\phi . \]
Let $D = (\phi^{-1} \nabla^h )^H$ i.e. $D V \in C^\infty (H(M)^* \otimes \phi^{-1} T N)$ is the restriction of $(\phi^{-1} \nabla^h )V$
to $H(M)$. An $L^2$ inner product on $C^\infty (H(M)^* \otimes \phi^{-1} T N)$ is given by
\[ (\varphi , \psi ) = \int_M \langle \varphi , \psi \rangle \; \Psi , \;\;\; \left. \langle \varphi , \psi \rangle \right|_U =
\sum_{a=1}^{2n} h^\phi (\varphi X_a , \psi X_a ) , \]
for any $\varphi , \psi \in C^\infty (H(M)^* \otimes \phi^{-1} T N)$ and any local orthonormal frame $\{ X_a : 1 \leq a \leq 2n \}$ of
$H(M)$ on $U \subseteq M$. Then for any $V \in C^\infty_0 (\phi^{-1} T N)$
\[ (D^* \varphi , V) = \int \sum_a h^\phi (\varphi X_a \, , \, (\phi^{-1} \nabla^h )_{X_a} V) \Psi = \]
\[ = \int \sum_a \{ X_a (h^\phi (\varphi X_a , V)) - h^\phi ((\phi^{-1} \nabla^h )_{X_a} \varphi X_a \, , \, V) \} \Psi . \]
Let $X_{\varphi , V} \in H(M)$ be determined by $G_\theta (X_{\varphi , V} , Y) = h^\phi (\varphi Y , V)$ for any $Y \in H(M)$. Then
\[ \sum_a X_a (h^\phi (\varphi X_a , V)) = \sum_a X_a (g_\theta (X_{\varphi , V} , X_a )) = \]
\[ = \sum_a \{ g_\theta (\nabla_{X_a} X_{\varphi , V} , X_a ) + g_\theta (X_{\varphi , V} , \nabla_{X_a} X_a )\} = \]
\[ = {\rm div}(X_{\varphi , V}) + \sum_a h^\phi (\varphi \nabla_{X_a} X_a \, , \, V) . \]
We may conclude that
\[ D^* \varphi = - \sum_{a=1}^{2n} \{ (\phi^{-1} \nabla^h )_{X_a} \varphi X_a - \varphi \nabla_{X_a} X_a \} \]
on $U$ and then $D^* D V = - \Delta_b^\phi V$ for any $V \in C^\infty (\phi^{-1} T N)$.
\begin{proposition}
The symbol of the rough sublaplacian is
\begin{equation}
\sigma_2 \left( \Delta_b^\phi \right)_\omega v = \left[ \omega (T_x )^2 - \| \omega \|^2 \right] \, v
\label{e:Symbol}
\end{equation}
for any $\omega \in T^\ast_x (M) \setminus \{ 0 \}$, $v \in \left( \phi^{-1} T N \right)_x$ and $x \in M$. Therefore $\Delta_b^\phi$ is
a degenerate elliptic operator and its ellipticity degenerates precisely in the cotangent directions spanned by $\theta$.
\end{proposition}
{\em Proof}. Let $T^\prime (M) = T^\ast (M) \setminus (0)$ and let $\Pi : T^\prime (M) \to M$ be the projection. If $E \to M$ and $F \to M$ are vector bundles
we set
\[ {\rm Smbl}_k (E, F) = \left\{ \sigma \in {\rm Hom}(\Pi^\ast E \, , \, \Pi^\ast F) : \sigma_{\rho \omega} = \rho^k \, \sigma_\omega \, , \; \rho > 0 \right\} \]
(with $k \in {\mathbb Z}$).
Let $\sigma_k (L) \in {\rm Smbl}_k (E, F)$ be the symbol of the $k$-th order differential operator $L \in {\rm Diff}_k (E,F)$. We wish to compute
$\sigma_2 (\Delta_b^\phi ) \in {\rm Smbl}_2 (\phi^{-1} T N \, , \, \phi^{-1} T N)$. To this end let $\omega \in T^\prime (M)$ such that $\Pi (\omega ) = x$
and let $f \in C^\infty (M)$ such that $(d f)_x = \omega$. Also let $v \in (\phi^{-1} T N)_x$ and $V \in C^\infty (\phi^{-1} T N)$ such that $V_x = v$. Then
\[ \sigma_2 (\Delta_b^\phi ) = - \frac{1}{2} \, \Delta_b^\phi \left[ (f - f(x))^2 V \right] (x). \]
Then (\ref{e:Symbol}) follows from the identities
\[ \Delta_b (u^2 ) = 2 u \, \Delta_b u + 2 \| \nabla^H u \|^2 \, , \]
\[ \Delta_b^\phi (g V) = g \, \Delta_b^\phi V + (\Delta_b g) V + 2 (\phi^{-1} \nabla^h )_{\nabla^H g} V \, , \]
where $g = u^2$ and $u = f - f(x)$. The norm in (\ref{e:Symbol}) is $\| \omega \| = g_{\theta , x}^\ast (\omega , \omega )$. Hence for each $v \in
{\rm Ker} \left[ \sigma_2 (\Delta_b^\phi )_\omega \right]$ either $v = 0$ or $\omega = \lambda \theta_x$ for some $\lambda \in {\mathbb R} \setminus \{ 0 \}$.
Q.e.d.

\section{The first variation formula}
Let $\tilde{M} = M \times (- \epsilon , \epsilon )$, $\epsilon > 0$, and let $F : \tilde{M} \to N$ be a smooth $1$-parameter
variation of $\phi$ by smooth maps i.e. $\phi_0 = \phi$ where $\phi_t = F \circ \alpha_t$ and $\alpha_t : M \to \tilde{M}$ is the injection
$\alpha_t (x) = (x,t)$ for any $x \in M$. Let $V \in C^\infty (\phi^{-1} TN)$ be the corresponding infinitesimal variation
\[ \left. V_x = (d_{(x,0)} F) \frac{\partial}{\partial t} \right|_{(x,0)} \, , \;\;\; x \in M. \]
Let us consider the {\em second fundamental form} of $\phi$ (cf. R. Petit, \cite{Pet})
\[ \beta_b (\phi) (X,Y) = \left( f^{-1} \nabla^h \right)_X \phi_* Y - \phi_* \nabla_X Y , \;\;\; X,Y \in {\mathfrak X}(M). \]
Here $\phi_* X$ denotes the cross-section in $\phi^{-1} T N \to M$ given by
\[ (\phi_* X) (x) = (d_x \phi) X_x \, , \;\;\; x \in M, \]
for each $X \in {\mathfrak X}(M)$. Then (cf. e.g. \cite{BaDrUr1})
\[ \tau_b (\phi ) = {\rm trace}_{G_\theta} \, \Pi_H \beta_b (\phi ) = \sum_{a=1}^{2n} \beta_b (\phi )(X_a , X_a ) \]
on $U \subset M$. We shall establish the following
\begin{theorem}
Let $M$ be a compact strictly pseudoconvex CR manifold and $\phi : M \to N$ a smooth map into a
Riemannian manifold. Then for any smooth $1$-parameter variation $F : \tilde{M} \to N$ of $\phi$
\begin{equation}
\frac{d}{dt} \left\{ E_{2,b} (\phi_t ) \right\}_{t=0} = \int_M h^\phi \left( V , \, BH_b (\phi )  \right) \; \Psi
\label{e:SubBih10}
\end{equation}
where $\hat{R}^h \equiv \phi^{-1} R^h$ is given by $\left( \hat{R}^h (u,v)w \right)_x = R^h_{\phi (x)} (u_x , v_x ) w_x$, for any $u,v,w
\in C^\infty (\phi^{-1} TN)$ and any $x \in N$.
\label{p:SubBih2}
\end{theorem}
If $M$ is not compact one may, as usual, integrate over an arbitrary re\-la\-ti\-vely compact domain $\Omega \subset M$ and
consider only smooth $1$-parameter variations supported in $\Omega$.
\par
{\em Proof of Theorem} \ref{p:SubBih2}. Given $X \in {\mathfrak X}(M)$ we define $\tilde{X} \in {\mathfrak X}(\tilde{M})$ by
setting
\[ \tilde{X}_{(x,t)} = (d_x \alpha_t )X_x , \;\;\; (x,t) \in \tilde{M} . \]
Then
\[ \left( (\phi_t )_* X \right)_x = (d_x \phi_t ) X_x = (d_{(x,t)} F)(d_x \alpha_t ) X_x = (d_{(x,t)} F) \tilde{X}_{(x,t)} \] hence
\begin{equation}
(\phi_t )_* X = (F_* \tilde{X}) \circ \alpha_t , \;\;\; |t| < \epsilon ,
\label{e:SubBih11}
\end{equation}
with the obvious meaning of $F_* \tilde{X}$ as a section in $F^{-1} T N \to \tilde{M}$. Similarly the identities
\begin{equation}
X_i^{\phi_t} = X_i^F \circ \alpha_t , \;\;\; 1 \leq i \leq \nu ,
\end{equation}
relate the local frames $\{ X_i^{\phi_t} : 1 \leq i \leq \nu \}$ and $\{ X_i^F : 1 \leq i \leq \nu \}$ in the pullback bundles
$\phi_t^{-1} T N \to M$ and $F^{-1} T N \to \tilde{M}$ respectively. Let $\{ X_a : 1 \leq a \leq 2 n \}$ be a local orthonormal frame of $H(M)$. We wish to compute
\[ \tau_b (\phi_t ) = \sum_{a=1}^{2n} \{ (\phi_t^{-1} \nabla^h )_{X_a} (\phi_t )_* X_a - (\phi_t )_* \nabla_{X_a} X_a \} \]
for any $|t| < \epsilon$. Note that
\begin{equation}
(\phi^{-1} \nabla^h )_X \phi_* Y = \{ X(Y \phi^i ) + X(\phi^j ) Y(\phi^k ) \left( \Gamma^i_{jk} \circ \phi \right) \} X_i^\phi
\label{e:SubBih13}
\end{equation}
where $\phi^i = y^i \circ \phi$. As $\phi_t^i = F^i \circ \alpha_t$ it follows that
\begin{equation} X(\phi_t^i ) = \tilde{X}(F^i ) \circ \alpha_t , \;\;\; X(Y \phi^i_t ) = [ \tilde{X} (\tilde{Y} F^i ) ] \circ \alpha_t ,
\label{e:SubBih14}
\end{equation}
for any $X,Y \in {\mathfrak X}(M)$. Therefore (by (\ref{e:SubBih13})-(\ref{e:SubBih14}))
\begin{equation}
\label{e:SubBih15} (\phi_t^{-1} \nabla^h )_X (\phi_t )_* Y = \left[ (F^{-1} \nabla^h )_{\tilde{X}} F_* \tilde{Y} \right] \circ \alpha_t
\end{equation}
where $F^{-1} \nabla^h$ is the connection in $F^{-1} T N \to \tilde{M}$ induced by $\nabla^h$. Then (by (\ref{e:SubBih11}) and (\ref{e:SubBih15}))
\begin{equation}
\tau_b (\phi_t ) = \tau_b (F) \circ \alpha_t , \;\;\; |t| < \epsilon ,
\label{e:SubBih16}
\end{equation}
where $\tau_b (F) \in C^\infty (F^{-1}(V), F^{-1} TN)$ is defined by
\[ \tau_b (F) = \sum_{a=1}^{2n} \{ (F^{-1} \nabla^h )_{\tilde{X}_a} F_* \tilde{X}_a  - F_* \widetilde{\nabla_{X_a} X_a} \} . \]
Let $h^F$ be the metric induced by $h$ in $F^{-1} TN \to \tilde{M}$ so that
\[ h^{\phi_t} (X_i^{\phi_t} \, , \, X_j^{\phi_t}) = h^F (X_i^F \, , \, X_j^F ) \circ \alpha_t , \;\;\; 1 \leq i,j \leq \nu . \]
Consequently
\[ \| \tau_b (\phi_t )\|^2 = \sum_{a=1}^{2n} h^F \left( (F^{-1} \nabla^h )_{\tilde{X}_a} F_* \tilde{X}_a - F_*
\widetilde{\nabla_{X_a} X_a} \, , \, \tau (F) \right) \circ \alpha_t \]
so that
\[ \frac{d}{dt} \left\{ E_{2,b} (\phi_t ) \right\}_{t=0} = \frac{1}{2} \int_M \frac{\partial}{\partial t} \left\{ h^F (\tau (F), \tau (F)) \right\}_{(x,0)} \; \Psi (x) = \]
\[ = \int_M \sum_{a=1}^{2n} h^F \left( (F^{-1} \nabla^h )_{\partial /\partial t} \left[ (F^{-1} \nabla^h )_{\tilde{X}_a}
F_* \tilde{X}_a - \right.\right. \]
\[ \left. \left. - F_* \widetilde{\nabla_{X_a} X_a} \right] \, , \, \tau_b (F) \right)_{(x,0)} \; \Psi (x) . \]
Let $R^{F^{-1} \nabla^h}$ be the curvature tensor field of $F^{-1} \nabla^h$. As
\[ \left[ \tilde{X}_a \, , \, \frac{\partial}{\partial t}\right] = 0, \;\;\; 1 \leq a \leq 2n, \]
it follows that
\begin{equation}
(F^{-1} \nabla^h )_{\partial /\partial t} (F^{-1} \nabla^h )_{\tilde{X}_a} F_* \tilde{X}_a =
\label{e:SubBih17}
\end{equation}
\[ = (F^{-1} \nabla^h )_{\tilde{X}_a} (F^{-1} \nabla^h )_{\partial /\partial t} F_* \tilde{X}_a - R^{F^{-1} \nabla^h} (\tilde{X}_a \,
, \, \frac{\partial}{\partial t}) F_* \tilde{X}_a . \]
On the other hand
\[ (F^{-1} \nabla^h )_{\partial /\partial t} F_* \tilde{X} - (F^{-1} \nabla^h )_{\tilde{X}} F_* \frac{\partial}{\partial t} = \]
\[ = \{ \frac{\partial}{\partial t} \left( \tilde{X} F^i \right) + \tilde{X}(F^j ) \frac{\partial F^k}{\partial t} \left(
\Gamma^i_{jk} \circ F \right) \} X_i^F - \]
\[ - \{ \tilde{X} \left( \frac{\partial F^i}{\partial t} \right) + \frac{\partial F^j}{\partial t} \tilde{X}(F^k ) \left( \Gamma^i_{jk} \circ F \right) \} X_i^F = \]
\[ = \left[ \frac{\partial}{\partial t} \, , \, \tilde{X}\right] (F^i ) X_i^F = F_* \left[ \frac{\partial}{\partial t} \, , \, \tilde{X} \right] = 0 \]
so that
\begin{equation} (F^{-1} \nabla^h )_{\partial /\partial t} F_* \tilde{X} = (F^{-1} \nabla^h )_{\tilde{X}} F_* \frac{\partial}{\partial t} .
\label{e:SubBih18}
\end{equation}
By (\ref{e:SubBih17})-(\ref{e:SubBih18})
\begin{equation}\label{e:SubBih19}
\frac{d}{dt} \left\{ E_{2,b} (\phi_t ) \right\}_{t=0} = \sum_{a=1}^{2n} \int_M h^F \left( (F^{-1}
\nabla^h )_{\tilde{X}_a} (F^{-1} \nabla^h )_{\tilde{X}_a} F_* \frac{\partial}{\partial t} - \right.
\end{equation}
\[ - \left. (F^{-1} \nabla^h )_{\widetilde{\nabla_{X_a} X_a}} F_* \frac{\partial}{\partial t} - R^{F^{-1} \nabla^h} (\tilde{X}_a ,
\frac{\partial}{\partial t}) F_* \tilde{X}_a \, , \, \tau_b (F) \right)_{(x,0)} \; \Psi (x) . \]
We compute separately the integrand in the right hand side of (\ref{e:SubBih19}) as follows. First (as $(F^{-1} \nabla^h ) h^F = 0$)
\[ h^F \left( (F^{-1} \nabla^h )_{\tilde{X}_a} (F^{-1} \nabla^h )_{\tilde{X}_a} F_* \frac{\partial}{\partial t} \, , \, \tau_b (F)\right) = \]
\[ = \tilde{X}_a \left( h^F \left( (F^{-1} \nabla^h )_{\tilde{X}_a} F_* \frac{\partial}{\partial t} \, , \, \tau_b (F) \right) \right) - \]
\[ - h^F \left( (F^{-1} \nabla^h )_{\tilde{X}_a} F_* \frac{\partial}{\partial t} \, , \, (F^{-1} \nabla^h )_{\tilde{X}_a} \tau_b (F) \right) = \]
\[ = \tilde{X}_a \tilde{X}_a \left( h^F \left( F_* \frac{\partial}{\partial t} \, , \, \tau_b (F) \right) \right) - \]
\[ - 2 \tilde{X}_a \left( h^F \left( F_* \frac{\partial}{\partial t} \, , \, (F^{-1} \nabla^h )_{\tilde{X}_a} \tau_b (F) \right) \right) + \]
\[ + h^F \left( F_* \frac{\partial}{\partial t} \, , \, (F^{-1} \nabla^h )_{\tilde{X}_a} (F^{-1} \nabla^h )_{\tilde{X}_a} \tau_b (F) \right) . \]
Moreover, as $\tilde{X}(\varphi ) \circ \alpha_t = X(\varphi \circ \alpha_t )$ for any $X \in {\mathcal X}(M)$ and any $\varphi \in C^\infty ( \tilde{M})$
\begin{equation}
\left[ \tilde{X}_a \tilde{X}_a \left( h^F \left( F_* \frac{\partial}{\partial t} \, , \, \tau (F) \right) \right) \right] \circ \alpha_t =
\label{e:SubBih20}
\end{equation}
\[ = X_a X_a \left( h^F \left( F_* \frac{\partial}{\partial t} \, , \, \tau_b (F) \right) \circ \alpha_t \right) . \]
Another tautology we shall need is
\[ V = \left( F_* \frac{\partial}{\partial t} \right) \circ \alpha_0 . \]
The following calculation
\[ \left[ (F^{-1} \nabla^h )_{\tilde{X}} \tau_b (F) \right]_{(x,t)} = \]
\[ = \{ \tilde{X}(\tau_b (F)^i ) + \tilde{X}(F^j ) \tau_b (F)^k \left( \Gamma^i_{jk} \circ F \right) \}_{(x,t)} X_i^F (x,t) = \]
\[ = \{ X(\tau_b (\phi_t )^i ) + X(\phi_t^j ) \tau_b (\phi_t )^k \left( \Gamma^i_{jk} \circ \phi_t \right) \}_x X_i^{\phi_t} (x) \]
shows that
\begin{equation}
\left[ (F^{-1} \nabla^h )_{\tilde{X}_a} \tau_b (F) \right] \circ \alpha_t = (\phi_t^{-1} \nabla^h )_{X_a} \tau_b (\phi_t ).
\label{e:SubBih21}
\end{equation}
Let $X_t \in H(M)$ be the horizontal tangent vector field determined by
\[ G_\theta (X_t , Y) = h^F \left( F_* \frac{\partial}{\partial t} \, , \, (F^{-1} \nabla^h )_{\tilde{Y}} \tau_b (F) \right) \circ \alpha_t \]
for any $Y \in H(M)$. Then (by $\nabla g_\theta = 0$)
\[ \tilde{X}_a \left( h^F \left( F_* \frac{\partial}{\partial t} \, , \, (F^{-1} \nabla^h )_{\tilde{X}_a} \tau_b (F) \right) \right)
\circ \alpha_t = X_a (G_\theta (X_t , X_a )) = \]
\[ = g_\theta (\nabla_{X_a} X_t , X_a ) + g_\theta (X_t \nabla_{X_a} X_a ) \]
that is
\begin{equation}
\tilde{X}_a \left( h^F \left( F_* \frac{\partial}{\partial t} \, , \, (F^{-1} \nabla^h )_{\tilde{X}_a} \tau_b (F) \right) \right) \circ
\alpha_t = g_\theta (\nabla_{X_a} X_t \, , \, X_a ) +
\label{e:SubBih22}
\end{equation}
\[ + h^F \left( F_* \frac{\partial}{\partial t} \, , \, (F^{-1} \nabla^h )_{\widetilde{\nabla_{X_a} X_a}} \tau_b (F) \right) \circ \alpha_t . \]
A calculation similar to that leading to (\ref{e:SubBih21}) furnishes
\begin{equation}
\left[ (F^{-1} \nabla^h )_{\tilde{X}} (F^{-1} \nabla^h )_{\tilde{Y}} \tau_b (F) \right] \circ \alpha_t = (\phi_t^{-1} \nabla^h
)_X (\phi_t^{-1} \nabla^h )_Y \tau_b (\phi_t ) .
\label{e:SubBih23}
\end{equation}
Summing up the information in (\ref{e:SubBih20}) and (\ref{e:SubBih22})-(\ref{e:SubBih23})
\begin{equation}
\label{e:SubBih24}
\sum_{a=1}^{2n} h^F \left( (F^{-1} \nabla^h )_{\tilde{X}_a} (F^{-1} \nabla^h )_{\tilde{X}_a} F_* \frac{\partial}{\partial t}
\, , \, \tau_b (F) \right) \circ \alpha_0 =
\end{equation}
\[ = h^\phi \left( V , \, \Delta_b^\phi \, \tau_b (\phi) \right) - 2 {\rm div}(X_0 ) + \]
\[ + \sum_a \{ X_a X_a \left( h^\phi (V, \, \tau_b (\phi ))\right) - h^\phi \left( V, \, (\phi^{-1} \nabla^h )_{\nabla_{X_a} X_a} \tau_b (\phi )\right) \} . \]
Using (\ref{e:SubBih24}) and
\[ h^F \left( (F^{-1} \nabla^h )_{\widetilde{\nabla_{X_a} X_a}} F_* \frac{\partial}{\partial t} \, , \, \tau_b (F) \right) = \]
\[ = \widetilde{\nabla_{X_a} X_a} \left( h^F \left( F_* \frac{\partial}{\partial t} \, , \, \tau_b (F) \right) \right) - h^F
\left( F_* \frac{\partial}{\partial t} \, , \, (F^{-1} \nabla^h )_{\widetilde{\nabla_{X_a} X_a}} \tau_b (F) \right) \]
we may conclude that
\begin{equation}
\sum_{a=1}^{2n} h^F \left( (F^{-1} \nabla^h )_{\tilde{X}_a} (F^{-1} \nabla^h )_{\tilde{X}_a} F_* \frac{\partial}{\partial t} -
\right. \label{e:SubBih25}
\end{equation}
\[ - \left. (F^{-1} \nabla^h )_{\widetilde{\nabla_{X_a} X_a}} F_* \frac{\partial}{\partial t} \, , \, \tau_b (F) \right) \circ \alpha_0 = \]
\[ = h^\phi \left( V , \, \Delta_b^\phi \, \tau_b (\phi ) \right) - 2 {\rm div}(X_0 ) +  \Delta_b \left[ h^\phi (V, \tau_b (\phi )) \right] . \]
Given local coordinates $(U, \tilde{x}^A )$ on $M$ we set $x^A  = \tilde{x}^A \circ p$ where $p : \tilde{M} \to M$ is the natural projection. To compute the
curvature term in the right hand side of (\ref{e:SubBih19}) we conduct
\[ \left( R^{F^{-1} \nabla^h} (\tilde{X} \, , \, \frac{\partial}{\partial t}) F_* \tilde{Y} \right)_{(x,0)} = \]
\[ = X^A (x) \tilde{Y}(F^j )_{(x,0)} \left( R^{F^{-1} \nabla^h} \left( \frac{\partial}{\partial x^A} \, , \,
\frac{\partial}{\partial t} \right) X^F_j \right)_{(x,0)} = \]
\[ = X^A (x) Y(\phi^j )_x \{ (F^{-1} \nabla^h )_{\partial /\partial x^A} (F^{-1} \nabla^h )_{\partial /\partial t} X_j^F - \]
\[ - (F^{-1} \nabla^h )_{\partial /\partial t} (F^{-1} \nabla^h )_{\partial /\partial x^A} X_j^F \}_{(x,0)} = \]
\[ = X^A (x) Y(\phi^j )_x \{ (F^{-1} \nabla^h )_{\partial /\partial x^A} [ \frac{\partial F^k}{\partial t} \left( \Gamma^i_{kj} \circ
F \right) X_i^F ] - \]
\[ - (F^{-1} \nabla^h )_{\partial /\partial t} [ \frac{\partial F^k}{\partial x^A} \left( \Gamma^i_{kj} \circ F \right) X_i^F ]
\}_{(x,0)} = \]
\[ = X^A (x) Y(\phi^j )_x \left\{ \left( \frac{\partial^2 F^k}{\partial x^A \partial t} (\Gamma^i_{kj}\circ F) +
\frac{\partial F^k}{\partial t} (\frac{\partial \Gamma^i_{kj}}{\partial y^\ell} \circ F) \frac{\partial
F^\ell}{\partial x^A} \right) X_i^F + \right. \]
\[ +  \frac{\partial F^k}{\partial t} (\Gamma^i_{kj}\circ F) \frac{\partial F^\ell}{\partial x^A} (\Gamma^m_{\ell i}\circ F) X_m^F -  \]
\[ - \left( \frac{\partial^2 F^k}{\partial t \partial x^A} (\Gamma^i_{kj}\circ F) + \frac{\partial F^k}{\partial x^A}
(\frac{\partial \Gamma^i_{kj}}{\partial y^\ell} \circ F) \frac{\partial F^\ell}{\partial t} \right) X_i^F - \] \[ - \left.
\frac{\partial F^k}{\partial x^A} (\Gamma^i_{kj}\circ F) \frac{\partial F^\ell}{\partial t} (\Gamma^m_{\ell i}\circ F) X_m^F \right\}_{(x,0)} = \]
\[ = X^A (x) Y(\phi^j )_x \frac{\partial F^k}{\partial t}(x,0) \frac{\partial F^\ell}{\partial x^A}(x,0) \times \]
\[ \times \left( \frac{\partial \Gamma^m_{kj}}{\partial y^\ell} - \frac{\partial \Gamma^m_{\ell j}}{\partial y^k} + \Gamma^i_{kj}
\Gamma^m_{\ell i} - \Gamma^i_{\ell j} \Gamma^m_{ki} \right)_{\phi (x)} X^\phi_m (x)  = \]
\[ = X^A (x) Y(\phi^j )_x \left( R^h \right)^m_{\ell kj} (\phi (x)) \frac{\partial F^k}{\partial t}(x,0) \frac{\partial
F^\ell}{\partial x^A}(x,0) X^\phi_m (x)  \]
where $X = X^A \partial /\partial \tilde{x}^A$ on $U$. Noticing that
\[ X(\phi^i )_x = X^A (x) \frac{\partial F^i}{\partial x^A}(x,0), \;\;\; x \in U, \]
we obtain
\begin{equation}
\left( R^{F^{-1} \nabla^h} (\tilde{X} \, , \, \frac{\partial}{\partial t}) F_* \tilde{Y} \right)_{(x,0)} =
\end{equation}
\[ = X(\phi^i )_x Y(\phi^j )_x \frac{\partial F^k}{\partial t}(x,0) \left( R^h \right)^m_{ikj} (\phi (x)) X_m^\phi (x). \]
Finally
\[ V_x = \frac{\partial F^i}{\partial t}(x,0) X_i^\phi (x), \;\;\ X \in U, \]
yields
\begin{equation}
\left( R^{F^{-1} \nabla^h} (\tilde{X} \, , \, \frac{\partial}{\partial t}) F_* \tilde{Y} \right)_{(x,0)} =
R^h_{\phi (x)} ((\phi_* X)_x , \, V_x ) (\phi_* Y)_x
\label{e:SubBih26}
\end{equation}
for any $x \in M$. Next (by (\ref{e:SubBih26}) and the symmetries of the Riemann-Christoffel 4-tensor of $(N, h )$)
\[ h^F \left( R^{F^{-1} \nabla^h} (\tilde{X}_a \, , \, \frac{\partial}{\partial t}) F_* \tilde{X}_a \, , \, \tau_b (F)\right)_{(x,0)} = \]
\[ = h_{\phi (x)} \left( R^h_{\phi (x)} ((\phi_* X_a )_x , \, V_x ) (\phi_* X_a )_x \, , \, \tau_b (\phi )_x \right) =  \]
\[ = - h_{\phi (x)} \left( R^h_{\phi (x)} (\tau_b (\phi )_x \, , \, (\phi_* X_a )_x  ) (\phi_* X_a )_x \, , \, V_x \right) . \]
Together with (\ref{e:SubBih25}) and Green's lemma this leads to the first variation formula (\ref{e:SubBih10}) for $E_{2,b}$ in Theorem \ref{p:SubBih2}.
\vskip 0.1in
The following concept is central to this paper. A smooth map $\phi : M \to N$ is said to be a {\em subeliptic biharmonic map} if $\phi$ is
a critical point of the functional $E_{2,b} : C^\infty (M, N) \to {\mathbb R}$. By Theorem \ref{p:SubBih2} a smooth map $\phi : M \to N$ is subelliptic biharmonic
if and only if $\phi$ is a solution to (\ref{e:SubBih4}). A pseudoharmonic map is trivially subelliptic biharmonic.

\section{Subelliptic biharmonic maps and Fefferman's metric} We start by proving the identity
${\mathbb E}_2 (\phi \circ \pi ) = 2 \pi E_{2,b} (\phi )$ in Theorem \ref{t:SubBih1}. Let us set $\partial_A = \lambda_A^B T_B$ with
$\lambda^B_A \in C^\infty (U)$. As to the range of indices, we adopt the convention
\[ A,B,C, \cdots \in \{ 0, 1, \cdots , n , \overline{1}, \cdots , \overline{n} \} \]
with $T_0 = T$. Then (by (\ref{e:SubBih1}))
\[ F_{AB} = \tilde{G}_\theta (\partial_A , \partial_B) + \theta (\partial_A) \sigma (\partial_B) + \theta (\partial_B) \sigma (\partial_A ) = \]
\[ = g_{\alpha\overline{\beta}} (\lambda_A^\alpha \lambda_B^{\overline{\beta}} + \lambda_B^\alpha
\lambda_A^{\overline{\beta}} ) + \lambda_A^0 \sigma_B + \lambda_B^0 \sigma_A \]
where $\sigma_A = \sigma (\partial_A )$. A calculation based on (\ref{e:SubBih2}) shows that
\[ \sigma_A = \frac{1}{n+2} \{ i \lambda^B_A ( \Gamma_{B\alpha}^\alpha - \frac{1}{2} g^{\alpha\overline{\beta}}
T_B (g_{\alpha\overline{\beta}} )) - \frac{\rho}{4(n+1)} \, \lambda^0_A \} \circ \pi \]
where $\Gamma^\beta_{B\alpha}$ are given by $\nabla_{T_B} T_\alpha = \Gamma_{B\alpha}^\beta T_\beta$. Moreover (by (\ref{e:SubBih1}))
\[ F_{A,2n+2} = 2 [(\pi^* \theta ) \odot \sigma ](\partial_A \, , \, \partial /\partial \gamma ) = \frac{1}{n+2} \, \lambda^0_A \, , \]
\[ F_{2n+2,2n+2} = 0. \]
Next, using $F^{ab} F_{bc} = \delta^a_c$ (with $a,b,c, \cdots \in \{ 1, \cdots , 2n+2 \}$) we find
\begin{equation}
\begin{cases} F^{AB} F_{BC} + \displaystyle{\frac{\lambda_C^0}{n+2} \, F^{A,2n+2}} = \delta^A_C
\cr F^{AB} \lambda^0_B = 0 \cr F^{2n+2,B} F_{BC} + \displaystyle{\frac{\lambda^0_C}{n+2} \, F^{2n+2,2n+2}} = 0 \cr
F^{2n+2,B} \lambda^0_B = n+2. \cr
\end{cases}
\label{e:SubBih27}
\end{equation}
Taking into account that $\partial \Phi^j /\partial \gamma = 0$ and $F^{AB} \lambda_B^0 = 0$ we have
\begin{equation}
F^{pq} \left( \Gamma^i_{jk} \circ \Phi \right) \frac{\partial \Phi^j}{\partial u^A} \frac{\partial \Phi^k}{\partial u^B} =
\label{e:SubBih28}
\end{equation}
\[ =  F^{AB} \left( \Gamma^i_{jk} \circ \Phi \right) \left\{ \lambda^\alpha_A \lambda^\beta_B T_\alpha (\phi^j ) T_\beta
(\phi^k ) + \lambda^{\alpha}_A \lambda^{\overline{\beta}}_B T_\alpha (\phi^j ) T_{\overline{\beta}}(\phi^k )  + \right. \]
\[ + \left. \lambda^{\overline{\alpha}}_A \lambda^\beta_B T_{\overline{\alpha}}(\phi^j ) T_\beta (\phi^k ) +
\lambda^{\overline{\alpha}}_A \lambda^{\overline{\beta}}_B T_{\overline{\alpha}}(\phi^j ) T_{\overline{\beta}}(\phi^k ) \right\} \circ \pi . \]
We need the following
\begin{lemma}
The {\rm (}reciprocal{\rm )} Fefferman metric is related to the {\rm (}reciprocal{\rm )} Levi form by
\begin{equation}
F^{AB} \lambda_A^\alpha \lambda^{\overline{\beta}}_B = g^{\alpha\overline{\beta}} ,
\label{e:SubBih29}
\end{equation}
\begin{equation}
F^{AB} \lambda_A^\alpha \lambda_B^\beta = 0.
\label{e:SubBih30}
\end{equation}
\label{l:SubBih2}
\end{lemma}
{\em Proof}. The identities (\ref{e:SubBih27}) may be written
\[ F^{AB} g_{\alpha\overline{\beta}} (\lambda_B^\alpha \lambda_C^{\overline{\beta}} + \lambda_C^\alpha \lambda^{\overline{\beta}}_B ) + F^{AB} \lambda^0_C \sigma_B +
\frac{1}{n+2} \, F^{A, 2n+2} \lambda^0_C = \delta^A_C , \]
\[ F^{AB} \lambda^0_B = 0, \]
\[ F^{2n+2,B} g_{\alpha\overline{\beta}} (\lambda_B^\alpha \lambda_C^{\overline{\beta}} + \lambda_C^\alpha \lambda^{\overline{\beta}}_B ) + (n+2) \sigma_C + \]
\[ + F^{2n+2, B} \lambda^0_C \sigma_B + \frac{1}{n+2} \, F^{2n+2,2n+2} \lambda^0_C = 0, \]
\[ F^{2n+2,B} \lambda^0_B = n+2. \]
If $\mu := \lambda^{-1}$ then (by the first of the previous four identities)
\[ \mu^A_D = (\frac{1}{n+2} \, F^{A,2n+2} + F^{AB} \sigma_B ) \delta^0_D + F^{AB} g_{\alpha\overline{\beta}} (\lambda^\alpha_B
\delta^{\beta + n}_D + \lambda^{\overline{\beta}}_B \delta^\alpha_D ) \]
yielding
\begin{equation}
\begin{cases} \mu^A_0 = \displaystyle{\frac{1}{n+2} \, F^{A, 2n+2} + F^{AB} \sigma_B} \cr \mu^A_\alpha = F^{AB} g_{\alpha\overline{\beta}}
\lambda^{\overline{\beta}}_B \cr \mu^A_{\beta + n} = F^{AB} g_{\alpha\overline{\beta}} \lambda^\alpha_B . \cr
\end{cases}
\label{e:SubBih31}
\end{equation}
The second and third of the identities (\ref{e:SubBih31}) lead to (\ref{e:SubBih29}) and (\ref{e:SubBih30}), respectively. Lemma
\ref{l:SubBih2} is proved.
\vskip 0.1in
By Lemma \ref{l:SubBih2} we may write (\ref{e:SubBih28}) as
\begin{equation}
\label{e:SubBih32}
F^{pq} \left( \Gamma^i_{jk} \circ \Phi \right) \frac{\partial \Phi^j}{\partial u^A} \frac{\partial \Phi^k}{\partial u^B} = 2
\left\{ \left( \Gamma^i_{jk} \circ \phi \right) g^{\alpha\overline{\beta}} T_\alpha (\phi^j )
T_{\overline{\beta}}(\phi^k ) \right\} \circ \pi .
\end{equation}
By a result of J.M. Lee (cf. \cite{Lee}, or Proposition 2.8 in \cite{DrTo}, p. 140)
\[ \square (u \circ \pi ) = \left( \Delta_b u\right) \circ \pi , \;\;\; u \in C^2 (M). \]
Hence (by
(\ref{e:SubBih32}))
\begin{equation}
\tau (\Phi ) = \tau_b (\phi ) \circ \pi .
\label{e:SubBih33}
\end{equation}
On the other hand if $\pi^{-1} (U) \approx U \times S^1$ is a local trivialization chart of the canonical circle bundle and $u
\in C^\infty (M)$ is a function supported in $U$ then
\[ \int_{C(M)} u \circ \pi \;\; d \, {\rm vol}(F_\theta ) = 2 \pi \int_M u \;\; \Psi \]
(by integration along the fibres of $S^1 \to C(M) \to M$, cf. e.g. (2.49) in \cite{DrTo}, p. 141). Hence (by a partition of unity argument)
\begin{equation}
\label{e:SubBih34} \int_{C(M)} \| \tau (\Phi )\|^2  \; d \, {\rm vol}(F_\theta ) = 2 \pi \int_M \| \tau_b (\phi )\|^2 \; \Psi .
\end{equation}
To prove the next statement in Theorem \ref{t:SubBih1} let $\{ \phi_t \}_{|t| < \epsilon}$ be a smooth $1$-parameter variation of $\phi$
($\phi_0 = \phi$) so that $\Phi_t = \phi_t \circ \pi$ is a $1$-parameter variation of $\Phi$. Therefore (by (\ref{e:SubBih34})) if $\Phi$
is biharmonic then $\phi$ is a critical point of $E_{2,b}$.
\par
The converse doesn't follow from (\ref{e:SubBih34}) but rather from the first variation formula for $E_{2,b}$. Indeed let $\phi : M \to N$ be a smooth solution to
(\ref{e:SubBih4}). A slight modification of Jiang Guoying's arguments (cf. \cite{JiGu}) leads to the following
\begin{lemma}
A smooth map $\Phi : C(M) \to N$ is biharmonic if and only if $\Phi$ is a solution to
\begin{equation}
BH(\Phi ) \equiv \square^\Phi \tau (\Phi ) + {\rm trace}_{F_\theta} \, \left\{ (\Phi^{-1} R^h )(\tau (\Phi ), \, \Phi_* \, \cdot ) \Phi_* \, \cdot \right\} = 0
\label{e:SubBih35}
\end{equation}
where
\begin{equation}
\square^\Phi {\mathfrak V} = \sum_{p=1}^{2n+2} \epsilon_p \{ (\Phi^{-1} \nabla^h )_{{\mathcal X}_p}^2 - (\Phi^{-1} \nabla^h
)_{\nabla^{C(M)}_{{\mathcal X}_p} {\mathcal X}_p} \} {\mathfrak V}
\end{equation}
is the rough Laplacian on $C(M)$. Here ${\mathfrak V} \in C^\infty (\Phi^{-1} T N)$ and $\{ {\mathfrak X}_p : 1 \leq p
\leq 2n+2 \}$ is a local $F_\theta$-orthonormal {\rm (}i.e. $F_\theta ({\mathfrak X}_p , {\mathfrak X}_q ) = \epsilon_p
\delta_{pq}$ with $\epsilon_1 = \cdots = \epsilon_{2n+1} = - \epsilon_{2n+2} = 1${\rm )} frame in $T(C(M))$. Also
$\nabla^{C(M)}$ is the Levi-Civita connection of $(C(M), F_\theta )$.
\end{lemma}
Let $\{ X_a : 1 \leq a \leq 2n \}$ be a local $G_\theta$-orthonormal frame in $H(M)$, defined on the open set $U
\subseteq M$. Then $\{ X_a^\uparrow , \; T^\uparrow \pm S : 1 \leq a \leq 2n \}$ is a local $F_\theta$-orthonormal frame of $T(C(M))$.
Here for any $X \in {\mathfrak X}(M)$ we denote by $X^\uparrow \in {\mathfrak X}(C(M))$ the horizontal lift of $X$ with respect to the connection $1$-form $\sigma$ on the
canonical circle bundle (thought of as a principal $S^1$-bundle over $M$). We recall that $X^\uparrow_z \in {\rm Ker}(\sigma )_z$
and $(d_z \pi ) X^\uparrow_z = X_{\pi (z)}$ for any $z \in C(M)$. Also
\[ S = \frac{n+2}{2} \; \frac{\partial}{\partial \gamma} \, . \]
Let $V \in C^\infty (\phi^{-1} T N)$ and ${\mathfrak V} = V \circ \pi$. Let $\Phi = \phi \circ \pi$. The identities
\[ (\Phi^{-1} \nabla^h )_{\partial /\partial \gamma} X_k^\Phi = 0, \]
\[ (\Phi^{-1} \nabla^h )_{\partial /\partial u^A} X_k^\Phi = [(\phi^{-1} \nabla^h )_{\partial /\partial x^A} X^\phi_k ] \circ \pi , \]
imply
\begin{equation}
(\Phi^{-1} \nabla^h )_{X^\uparrow} {\mathfrak V} = [ (\phi^{-1} \nabla^h )_X V ] \circ \pi , \label{e:SubBih37}
\end{equation}
\begin{equation}
(\Phi^{-1} \nabla^h )_{X^\uparrow} (\Phi^{-1} \nabla^h )_{Y^\uparrow} {\mathfrak V} = [ (\phi^{-1} \nabla^h )_X (\phi^{-1} \nabla^h )_Y V ] \circ \pi ,
\label{e:SubBih38}
\end{equation}
for any $X,Y \in T(M)$. At this point we need
\begin{lemma}
For any $X,Y \in H(M)$
\[ \nabla^{C(M)}_{X^\uparrow} Y^\uparrow = (\nabla_X Y )^\uparrow - (d \theta )(X,Y) T^\uparrow - [A(X,Y) + (d \sigma )(X^\uparrow , Y^\uparrow )] S , \]
\[ \nabla^{C(M)}_{X^\uparrow} T^\uparrow = (\tau X + q X)^\uparrow , \]
\[ \nabla^{C(M)}_{T^\uparrow} X^\uparrow = (\nabla_T X + q X)^\uparrow + 2(d \sigma )(X^\uparrow , T^\uparrow ) S, \]
\[ \nabla^{C(M)}_{X^\uparrow} S = \nabla^{C(M)}_S X^\uparrow = (J X)^\uparrow , \]
\[ \nabla^{C(M)}_{T^\uparrow} T^\uparrow = Q^\uparrow , \;\;\; \nabla^{C(M)}_S S = 0, \]
\[ \nabla^{C(M)}_S T^\uparrow = \nabla^{C(M)}_{T^\uparrow} S  = 0, \]
where $q : H(M) \to H(M)$ is given by $G_\theta (q X , Y) =
(d \sigma )(X^\uparrow , Y^\uparrow )$ and $Q \in H(M)$ is given
by $G_\theta (Q , Y) = 2(d \sigma )(T^\uparrow , Y^\uparrow )$.
Also $\tau$ is the pseudohermitian torsion of $\nabla$ and $A(X,Y)
= G_\theta (\tau X , Y)$. \label{l:SubBih4}
\end{lemma}
Cf. Lemma 2 in E. Baletta et al., \cite{BaDrUr2}, p. 083504-26. As
a consequence of Lemma \ref{l:SubBih4}
\[ \nabla^{C(M)}_{X^\uparrow} X^\uparrow = (\nabla_X X )^\uparrow - A(X,X) S, \]
\[ \nabla^{C(M)}_{T^\uparrow \pm S} (T^\uparrow \pm S) = Q^\uparrow , \]
hence (by (\ref{e:SubBih37}))
\[ \sum_a (\Phi^{-1} \nabla^h )_{\nabla^{C(M)}_{X_a^\uparrow} X_a^\uparrow} {\mathfrak V} = \] \[ = \sum_a \{ (\Phi^{-1} \nabla^h
)_{(\nabla_{X_a} X_a )^\uparrow} - A(X_a , X_a )(\Phi^{-1} \nabla^h )_S \} {\mathfrak V}  = \]
\[ = \sum_a [(\phi^{-1} \nabla^h )_{\nabla_{X_a} X_a} V ] \circ \pi \]
as ${\rm trace}_{G_\theta} A = {\rm trace}(\tau ) = 0$ (cf. e.g. (1.59) in \cite{DrTo}, p. 37). Together with (\ref{e:SubBih38})
and
\[ (\Phi^{-1} \nabla^h )_S {\mathfrak V} = 0 \]
this allows one to conduct the following calculation
\[ \square^\Phi {\mathfrak V} = \sum_{a=1}^{2n} \{ (\Phi^{-1} \nabla^h )^2_{X_a^\uparrow} - (\Phi^{-1} \nabla^h
)_{\nabla^{C(M)}_{X_a^\uparrow} X_a^\uparrow} \} {\mathfrak V} + \]
\[ + \{ (\Phi^{-1} \nabla^h )^2_{T^\uparrow + S} - (\Phi^{-1} \nabla^h )_{\nabla^{C(M)}_{T^\uparrow + S} T^\uparrow + S} \}
{\mathfrak V} - \]
\[ - \{ (\Phi^{-1} \nabla^h )^2_{T^\uparrow - S} - (\Phi^{-1} \nabla^h )_{\nabla^{C(M)}_{T^\uparrow - S} T^\uparrow - S} \}
{\mathfrak V} = \]
\[ = \sum_a \{ (\phi^{-1} \nabla^h )^2_{X_a} V - (\phi^{-1} \nabla^h )_{\nabla_{X_a} X_a} V \} \circ \pi  \]
that is
\begin{equation}
\square^\Phi {\mathfrak V} = (\Delta_b^\phi \, V) \circ \pi .
\label{e:SubBih39}
\end{equation}
It remains that we compute the curvature term in (\ref{e:SubBih35}). As $\Phi_* X^\uparrow = (\phi_* X) \circ \pi$ for
any $X,Y \in {\mathcal X}(M)$
\[  (\Phi^{-1} R^h )(\tau (\Phi ) , \Phi_* X^\uparrow ) \Phi_* Y^\uparrow  = [ (\phi^{-1} R^h )(\tau_b (\phi ) , \phi_* X ) \phi_* Y ] \circ \pi \]
hence
\[ {\rm trace}_{F_\theta} \, \{ (\Phi^{-1} R^h )(\tau (\Phi ), \Phi_* \, \cdot ) \Phi_* \, \cdot \} = \]
\[ = \sum_a (\Phi^{-1} R^h )(\tau (\Phi ) , \Phi_* X_a^\uparrow ) \Phi_* X_a^\uparrow + \]
\[ + (\Phi^{-1} R^h )(\tau (\Phi ) , \Phi_* (T^\uparrow + S)) \Phi_* (T^\uparrow + S) - \] \[ - (\Phi^{-1} R^h )(\tau (\Phi ) ,
\Phi_* (T^\uparrow - S)) \Phi_* (T^\uparrow - S) = \]
\[ = {\rm trace}_{G_\theta}\, \pi_H \{ (\phi^{-1} R^h )(\tau_b (\phi ), \phi_* \, \cdot ) \phi_* \, \cdot \}  + \]
\[ + 2 \{ (\Phi^{-1} R^h )(\tau (\Phi ) , \Phi_* T^\uparrow )\Phi_* S + (\Phi^{-1} R^h )(\tau (\Phi ) , \Phi_* S) \Phi_* T^\uparrow \} \]
and $\Phi_* S = 0$ so that
\begin{equation}
{\rm trace}_{F_\theta} \, \{ (\Phi^{-1} R^h )(\tau (\Phi ), \Phi_* \, \cdot ) \Phi_* \, \cdot \} =
\label{e:SubBih40}
\end{equation}
\[ = {\rm trace}_{G_\theta}\, \pi_H \{ (\phi^{-1} R^h )(\tau_b (\phi ), \phi_* \, \cdot ) \phi_* \, \cdot \} . \]
Finally (by (\ref{e:SubBih39})-(\ref{e:SubBih40})) $\phi \circ \pi$ is biharmonic. Theorem \ref{t:SubBih1} is proved.

\section{Conclusions and open problems}
We introduced the new concept of a subelliptic biharmonic map as a $C^\infty$ solution to the system (\ref{e:SubBih4}). This is a quasilinear
system of variational origin whose principal part is the bi-sublaplacian $\Delta_b^2$. $\Delta_b^2$ is a fourth order hypoelliptic operator, though not
elliptic, so that our work is part of the program outlined in \cite{JoXu}. Higher order degenerate elliptic equations, say of the form $\Delta_b^k u = 0$, were not
studied in the present day PDEs literature (cf. e.g. \cite{ArCrLi} for the elliptic case). Our main result is a geometric interpretation of subelliptic
biharmonic maps within Lorentzian geometry i.e. each $C^\infty$ solution to (\ref{e:SubBih4}) may be characterized as the base map corresponding to a
$S^1$-invariant biharmonic map from the total space $C(M)$ of the canonical circle bundle endowed with the Fefferman metric. Although, as shown in $\S \, 4$,
it makes sense to look for weak solutions to (\ref{e:SubBih4}) our methods in this paper are purely geometric and
a study of local properties of weak subelliptic biharmonic maps appears nowhere in the mathematical literature. Neither may the partial regularity theory
be naively reduced to that of $S^1$-invariant biharmonic maps, as $C(M)$ is Lorentzian so that no natural distance function on $C(M)$ is available {\em a priori}.
If $\Omega \subset M$ is a bounded domain the functionals $E_{2, b}(\phi )$ and $\int_\Omega |\Delta_b \phi |^2 \; \Psi$ (as introduced by
S-Y.A. Chang \& L. Wang \& P.C. Yang, \cite{ChWaYa}, in the elliptic case for maps $\phi : \Omega \to S^\nu$) have not been compared so far (and of course the
corresponding regularity for $\delta \, \int_\Omega |\Delta_b \phi |^2 \; \Psi = 0$ is unknown).

\end{document}